\documentclass[reqno,12pt,twoside]{article}
\usepackage{fullpage}
\usepackage{amssymb}
\usepackage{amsmath}
\usepackage{amsthm}
\newtheorem{theorem}{Theorem}
\newtheorem{problem}{Problem}

\newtheorem{remark}{Remark}

\begin{document}

\title{\Large On a problem of  Nathanson}\vskip -1cm
\author{\large Yong-Gao Chen$^1$ and Min Tang$^2$}\date{} \maketitle
\vskip -4cm
\begin{center}\vskip -1cm
{\small 1. School of Mathematical Sciences and Institute of Mathematics,\\
Nanjing Normal University,
Nanjing 210023, P. R. CHINA\\
  ygchen@njnu.edu.cn }
 \end{center}

 \begin{center}
{ \small 2. School of Mathematics and Computer Science,\\ Anhui Normal University,
Wuhu 241003, P. R. CHINA\\
tmzzz2000@163.com\
}
\end{center}

\begin{abstract} A set $A$ of nonnegative integers is an asymptotic basis
of order $h$ if every sufficiently large integer can be
represented as the sum of $h$ integers (not necessarily distinct)
of $A$. An asymptotic basis $A$ of order $h$ is minimal if no
proper subset of $A$ is an asymptotic basis of order $h$. In this
paper, we resolve a problem of Nathanson on minimal asymptotic
bases of order $h$.

{\bf Keywords:} minimal asymptotic basis; partition; Nathanson's
problem; binary expansion

{\bf 2010 Mathematics Subject Classification:}  11B13
\end{abstract}

\section{Introduction}

Let $\mathbb{N}$ denote the set of all nonnegative integers and
$h$ be an integer with $h\geq2$. For $A \subseteq \mathbb{N}$ and
$n\in \mathbb{N}$, let
$$r_h(A,n)=\sharp \{ (a_1, a_2,\ldots,a_h)\in A^{h}: a_1+a_2+\cdots+a_h=n\}. $$
A set $A$ is called \emph{an asymptotic basis of order $h$} if
$r_h(A,n)\geq 1$
 for all sufficiently large integers $n$.
In 1955, St\"{o}hr \cite{Stohr} introduced the concept of minimal
asymptotic basis. An asymptotic basis $A$ of order $h$ is
\emph{minimal} if no proper subset of $A$ is an asymptotic basis
of order $h$. This means that, for any $a\in A$, the set
$E_{a}=hA\setminus h(A\setminus \{a\})$ is infinite.

In 1956, H\"{a}rtter \cite{Hartter} proved that minimal asymptotic
bases of order $h$ exist for all $h\geq2$. Nathanson
\cite{Nathason1974} constructed explicitly a minimal asymptotic
basis of order $2$ and an asymptotic basis of order $2$ no subsets
of which is minimal. Chen and Chen \cite{Chen} resolved some
problems on minimal asymptotic bases asked by Nathanson
\cite{Nathason1988}.  Jia and Nathanson \cite{Jia1989} gave an
explicit construction of minimal asymptotic bases of order $h$ for
every $h\geq 2$. For related problems concerning minimal
asymptotic bases, one may see
\cite{ErdosNathanson1978Proc}-\cite{ErdosNathanson1988Illinois},
\cite{JanczakSchoen2010JNT}, \cite{Lee}-\cite{Tang} and
\cite{NathansonSarkozy1989Proc}.

For any  nonempty subset $W$ of $\mathbb{N}$, denote by
$\mathcal{F}^{\ast}(W)$ the set of all finite, nonempty subsets of
$W$. Let $A(W)$ be the set of
 all numbers of the form $\sum\limits_{f\in F}2^{f}$, where $F\in \mathcal{F}^{\ast}(W)$.

In 1988, Nathanson \cite{Nathason1988} posed the following problem
(see also Jia and Nathanson \cite{Jia1989}).

\begin{problem}\label{prob1} If $\mathbb{N}=W_0\cup\cdots\cup W_{h-1}$ is a partition such that
$w\in W_i$ implies either $w-1\in W_i$ or $w+1\in W_i$, then is
$$A=A(W_0)\cup\cdots\cup A(W_{h-1})$$ a minimal asymptotic basis of
order $h$?
\end{problem}

 In 2011, Chen and Chen
\cite{Chen} obtained the following result.

\noindent{\bf Theorem A} {\it Let $h\geq 2$ and $t$ be the least
integer with $t>\log h/\log 2$. Let $\mathbb{N}=W_0\cup\cdots\cup
W_{h-1}$ be a partition such that each set $W_i$ is infinite and
contains $t$ consecutive integers for $i=1,\ldots,h$. Then
$$A=A(W_0)\cup\cdots\cup A(W_{h-1})$$ is  a minimal asymptotic basis
of order $h$.}

 By Theorem A,  the answer to
Problem \ref{prob1} is affirmative for $2\le h<4$. For $a<b$, let
$[a, b]$ denote the set of all integers in the interval $[a, b]$.
In this paper, the following result is proved. Thus the answer to
Problem \ref{prob1} is negative for $h\ge 4$.

\begin{theorem}\label{mainthm1} Let $h$ and $t$ be integers with $2\leq t\le \log h/\log 2$.
Then there exists a partition $\mathbb{N}=W_0\cup\cdots\cup
W_{h-1}$ such that each set $W_i$ is a union of  infinitely many
intervals of at least $t$ consecutive integers and
$$A=A(W_0)\cup\cdots\cup A(W_{h-1})$$ is not a minimal asymptotic basis of order $h$.
\end{theorem}

\begin{remark} For $2\leq t<\log h/\log 2$, the following stronger result is proved:
 there exists a partition $\mathbb{N}=W_0\cup\cdots\cup W_{h-1}$
such that each set $W_i$ is a union of  infinitely many intervals
of at least $t$ consecutive integers and $n\in hA(W_0)$ for all
sufficiently large integers $n$.
\end{remark}

\section{Proof of the theorem }

 Since $t\ge 2$, it follows that $h\ge 2^t\ge 4$. For any subset $X$ of $\mathbb{N}$,
 let $2^X=\{ 2^x : x\in X\}$.  Let $\{ m_i \}_{i=1}^\infty$ be a sequence of
integers with $m_1>2^{h+4}$ and $m_{i+1}-m_i>2^{h+4}$ $(i\ge 1)$.
Let
$$W_0=[0, m_1] \cup \left( \bigcup_{i=1}^\infty [m_i+t+1, m_{i+1}]  \right) $$
and
$$W_j =\bigcup_{\substack{i=1\\i\equiv j\hskip -2mm\pmod{h-1}}}^\infty [m_i+1, m_{i}+t],
\quad j=1,2,\dots , h-1.$$ It is clear that $$\mathbb{N}=W_0\cup
W_1\cup \cdots \cup W_{h-1}.$$ If $w\in \mathbb{N}\setminus W_0$,
then, by the definition of $W_i$, we have $w>m_1>2^{h+4}$ and
$w-t\in W_0$. Write
$$A=A(W_0)\cup\cdots\cup A(W_{h-1}).$$
For any positive integer, let the binary expansion of $n$ be
\begin{equation}\label{eq1} n=\sum_{f\in F_n} 2^f.\end{equation}

Divide into two cases according to $h>2^t$ and $h=2^t$.

{\bf Case 1:} $h>2^t$.

In this case, we will prove that all integers $n$ with $n\ge
h2^{h(2t+1)}$ are in $hA(W_0)$. Thus $A$ is not a minimal
asymptotic basis of order $h$.

Let $n\ge h2^{h(2t+1)}$. Now we split terms in the summation of
\eqref{eq1}. First, we split all $2^f$ with $f\in F_n\setminus
W_0$ into $2^t$ terms $2^{f-t}$. Then all terms are in $2^{W_0}$
and each term repeats at most $2^t+1$ times in the summation. We
continue to split terms in the summation. For any term $2^w$ in
the summation,  if $w>2t+1$ and none of $2^{w-i}$ $(1\le i\le
2t+1)$ appears in the summation, we split $2^w$ (split one of
$2^w$ if there are several such terms) as follows:

(a) $2^w=2^{w-1}+2^{w-1}$ if $w-1\in W_0$;

(b) $2^w=(2^t+1) 2^{w-t-1}+\cdots +(2^t+1) 2^{w-2t+1} + (2^{t}+1)
2^{w-2t}+2 \times 2^{w-2t-1}$ if $w-1\notin W_0$.

In Case (b), by the definition of $W_0$ and $w\in W_0$, we know
that the integers $w-t-i$ $(1\le i\le t+1)$ are all in $W_0$.

Since each split increases the number of terms at least one, the
splitting  procedure must be terminated in finite steps.   In the
final summation,  all terms are in $2^{W_0}$ and each term repeats
at most $2^t+1$ times. If $2^w$ appears, then  at least one of
$2^{w-i}$ $(1\le i\le 2t+1)$ appears. Let the final summation be
$$n=\sum_{j=1}^s 2^{w_s}$$
with $0\le w_1\le w_2\le \cdots \le w_s$.  Let $w_0=0$. Thus
$$0\le w_{i+1}-w_i\le w_{i+1}-(w_{i+1}-2t-1)=2t+1,\quad i=0, 1, \dots , s-1.$$
Since
$$h2^{h(2t+1)}\le n=\sum_{j=1}^s 2^{w_s}\le (2^t+1)\sum_{w=0}^{w_s} 2^w =(2^t+1)(2^{w_s+1}-1)<h2^{w_s+1},$$
it follows that $w_s\ge h(2t+1)$. On the other hand,
$$w_s=\sum_{i=0}^{s-1} (w_{i+1}-w_i)\le s(2t+1).$$
Hence $s\ge h$. Noting that $2^t+1\le h$ and $s\ge h$,  we can
split the final summation into $h$ nonempty sums such that all
terms in each nonempty sum are distinct. So $n\in hA(W_0)$.

{\bf Case 2:} $h=2^t$.

It is clear that $4\in A(W_0)$. Now we prove that $E_4=hA\setminus
{h(A\setminus \{4\})}$ is a finite set. Thus $A$ is not a minimal
asymptotic basis of order $h$.

Let $n>m_2$. We will show that $n\in h(A\setminus \{4\})$. That
is, $n\notin E_4$. Divide into the following subcases:

 {\bf Subcase 2.1:} $F_n\cap W_0\not= \{ 2\} $.

{\bf Subcase 2.1.1:} $F_n\cap W_0\not=\emptyset $ and $
|F_n\setminus W_0|\ge h-1$. Then $F_n\setminus W_0$ has a
partition
$$F_n\setminus W_0=L_1\cup L_2\cup \cdots \cup L_{h-1},$$ where
$L_i\not= \emptyset$ $(1\le i\le h-1)$ and for every $L_i$ there
exists a $W_j$ $(j\ge 1)$ with $L_i\subseteq W_j$. Let
$L_0=F_n\cap W_0$ and
$$a_i=\sum_{l\in L_i} 2^l, \quad 0\le i\le h-1.$$
Then $$a_i\in A\setminus \{ 4\},\quad 0\le i\le h-1 $$ and
$$n=a_0+\cdots +a_{h-1}.$$ Hence $n\in h(A\setminus \{4\})$.

{\bf Subcase 2.1.2:} $F_n\cap W_0\not=\emptyset $ and $ 1\le
|F_n\setminus W_0|\le h-2$. Let
$$F_n\setminus W_0 =\{ f_0, \dots , f_{l-1}\}$$
with $f_0>\cdots >f_{l-1}$. Then $f_0\ge m_1+1>2^{h+4}$. Let
$$ f_i= f_0-(i-l+1), \quad l\le i\le h-2$$ and $f_{h-1}=f_{h-2}$.
Put
$$a_0=\sum_{f\in F_n\cap W_0} 2^f, \quad a_i=2^{f_i},\quad 1\le i\le
h-1.$$ Since $$f_l> f_{l+1}> \cdots > f_{h-2}=f_{h-1}>2^{h+4}
-(h-2-l+1)\ge 2^{h+4} -(h-2)>2,$$  it follows that $$a_i\in
A\setminus \{ 4\},\quad 0\le i\le h-1
$$ and
$$n=a_0+\cdots +a_{h-1}.$$ Hence $n\in h(A\setminus \{4\})$.

{\bf Subcase 2.1.3:} $F_n\cap W_0\not=\emptyset $ and
$F_n\setminus W_0=\emptyset$. That is, $F_n\subseteq W_0$. Let
$$F_n=\{ g_0, \dots , g_{k-1}\}$$
with $g_0>\cdots >g_{k-1}$. Since $$n>m_2>2^{h+4}>1+2+2^2+\cdots
+2^{h+3},$$ we have $g_0\ge h+4$.

If $k=1$, then $F_n=\{ g_0\}$. Let $$a_i=2^{g_0-i-1},\quad 0\le
i\le h-2$$ and $a_{h-1}=a_{h-2}$. Since $$a_0> a_1> \cdots
>a_{h-2}= a_{h-1}=2^{g_0-h+1}> 4,$$ it follows that $$a_i\in A\setminus \{
4\},\quad 0\le i\le h-1 $$ and $$n=a_0+\cdots +a_{h-1}.$$ Hence
$n\in h(A\setminus \{4\})$.

If $k\ge 2$ and $g_1>2$, then we take
$$a_0=2^{g_1}+\cdots +2^{g_{k-1}},$$
$$a_i=2^{g_0-i}, \quad 1\le i\le h-2$$ and $a_{h-1}=a_{h-2}$.
Since $$a_1> \cdots > a_{h-2}=a_{h-1}=2^{g_0-h+2}>4,$$ it follows
that
$$a_i\in A\setminus \{
4\},\quad 0\le i\le h-1 $$ and $$n=a_0+\cdots +a_{h-1}.$$ Hence
$n\in h(A\setminus \{4\})$.

If $k\ge 2$ and $g_1=2$, then we take
$$a_0=2^{g_0-h+1}+2^{g_1}+\cdots +2^{g_{k-1}}$$
and $$a_i=2^{g_0-i},\quad 1\le i\le h-1.$$ Since $$a_1> \cdots
> a_{h-1}=2^{g_0-h+1}>4,$$ it follows that $$a_i\in A\setminus \{
4\},\quad 0\le i\le h-1 $$ and $$n=a_0+\cdots +a_{h-1}.$$ Hence
$n\in h(A\setminus \{4\})$.

{\bf Subcase 2.2:} $F_n\cap W_0= \{ 2\} $. By $n>m_2$, we have
$F_n\setminus W_0\not=\emptyset$. If $f\in F\setminus W_0$, then
$f>m_1>2^{h+4}$ and $f-t\in W_0$ (see the arguments before Case
1). Let
$$ a_0=2^2+\sum_{f\in F_n\setminus \{ 2\} } 2^{f-t}$$ and
$$a_1=\cdots =a_{h-1}=\sum_{f\in F_n\setminus \{ 2\} } 2^{f-t}.$$
Then $$a_i\in A(W_0)\setminus \{ 4\},\quad 0\le i\le h-1$$ and
$$n=a_0+\cdots +a_{h-1}$$ by $h=2^t$. Hence $$n\in h(A(W_0)\setminus
\{4\}).$$

This completes the proof of Theorem \ref{mainthm1}.

\medspace

\noindent\textbf{Acknowledgments.} The authors are supported by
the National Natural Science Foundation of China, Nos.11771211 and
11471017.  The first author is also supported by the Priority
Academic Program Development of Jiangsu Higher Education
Institutions.

\end{document}